\documentclass[12pt]{amsart}
\usepackage{amsmath,amssymb,amsthm,amsfonts}
\usepackage{mathrsfs}        
\usepackage{comment}
\usepackage{pifont}          
\usepackage{enumerate}
\usepackage{tikz}
\usepackage[font=small,labelfont={sf,bf}]{caption}

\newtheorem{tm}{Theorem}[section]
\newtheorem{lemma}[tm]{Lemma}

\newtheorem{cor}[tm]{Corollary}

\theoremstyle{definition}

\newtheorem{defi}[tm]{Definition}

\theoremstyle{remark}
\newtheorem{rem}[tm]{Remark}

\newcommand{\R}{\mathbb R}         
\newcommand{\Rd}{{\mathbb R}^d}

\newcommand{\C}{\mathbb C}         

\newcommand{\W}{\dot{W}}           
\newcommand{\Hs}{\dot{H}}

\def\al{\alpha}                    
\def\be{\beta}

\def\De{\Delta}
\def\la{\lambda}

\def\th{\theta}

\def\eps{\epsilon}

 \def\cS{\mathcal{S}}
 \def\cD{\mathcal{D}}

 \def\F{\mathscr{F}}

\def\l{\langle}
\def\r{\rangle}

\def\wh{\widehat}
\def\t{\tilde}

\title{Strichartz estimates for the vibrating plate equation}
\author{Elena Cordero and Davide Zucco}
\address{Department of Mathematics,  University of Torino,
Via Carlo Alberto 10, 10123
Torino, Italy}
\email{elena.cordero@unito.it}
\email{davide.zucco@unito.it}


\date{}


\begin{document}

\begin{abstract}
We study the dispersive properties of the linear vibrating plate
(LVP) equation. Splitting it into two Schr\"odinger-type equations
we show its close  relation with the Schr\"odinger equation. Then,
the homogeneous Sobolev spaces appear to be the natural  setting
to show
 Strichartz-type estimates for the LVP equation.

By showing a  Kato-Ponce inequality for homogeneous Sobolev spaces
we prove the well-posedness of the Cauchy problem for the LVP
equation with  time-dependent potentials.  Finally, we exhibit the
sharpness of our results. This is achieved by finding a suitable
solution for the stationary homogeneous vibrating plate equation.
\end{abstract}

\maketitle

\section{Introduction}
In this paper we consider the Cauchy problem of the \emph{vibrating plate equation}
\begin{equation}\label{cpv}
\begin{cases}
\; \partial^2_t u + \Delta^2 u=F(t,x)\\
\;u(0,x)=u_0(x),\,\,
\partial_t u (0,x)=u_1(x)
\end{cases}
\end{equation}
where $t\in \R,\,x\in \Rd$, $u,\,F:\, \R\times\R^d \rightarrow \C$, $\De =\sum_{j=1}^d \partial_{x_i}^2$ is the classical Laplace
operator and $\De^2 u = \De(\De u)$. This equation is also known under the name of
\emph{Germain-Lagrange equation} by the ones that first discovered
the correct differential equation as a model for the vibration of
an elastic surface. It is at the basis of
 the theory of elasticity and has applications in architecture and engineering, see for example \cite{KV}.
Related works in the framework of variational calculus can be found in the books of Gelfand and Fomin \cite{GeFo}, Gould \cite{gould} and Weinstock \cite{wein}.
Recent studies in Sobolev, Gevrey and modulation spaces by means of pseudo-differential operators and time-frequency techniques are contained  in
 \cite{AgliardiCicognani,agza,pau} and \cite{corzuc}.
 The $2$-dimensional study of the vibrations of a nonlinear elastic plate (Von K\'arm\'an equations) was developed in \cite{tataru}.
 For related problems see also \cite{LL91}.

In this paper we focus on dispersive properties of the vibrating plate equation. The study of  dispersive properties of evolution equations  have become of great importance in PDE,   with applications to local and global existence for nonlinear analysis, well-posedness in Sobolev spaces of lower order, scattering theory and many other topics (see, e.g., Tao's book \cite{tao} and the references therein).
The matter of fact is provided by the celebrated Strichartz estimates (see Section \ref{ST} below), that play the  principal role in the whole study.
\vskip0.1truecm

The \emph{free} vibrating plate equation
\begin{equation}\label{cpvom}
  \partial_t^2 u(t,x) + \De^2 u(t,x) = 0
\end{equation}
can be factorized as the following product
\[
 (\partial_t^2 + \De^2 )u = (i\partial_t + \De) (-i\partial_t + \De) u
\]
which displays  the interesting relation with the Schr\"odinger
equation. Namely, the vibrating plate operator $P=\partial_t^2 +
\De^2$ can be recovered by  composing two Schr\"odinger-type
operators: $S_1=i\partial_t + \De$ and $S_2=-i\partial_t + \De$.
This fact suggests to recover Strichartz estimates for the
dispersive problem \eqref{cpv} from the  well-known ones  for the
Schr\"odinger equation, recalled shortly in Section \ref{secsch}.
Classical references on the subject are provided by \cite{gvelo}
and \cite{ktao}, see also \cite{cazenave}. The $L^p$ environment
of the estimates for the Schr\"odinger equation yields to
homogeneous Sobolev spaces as natural setting for the study of
\eqref{cpv}.

We shall use the Strichartz estimates for the VP equation to prove
a well-posedness result for the Cauchy problem \eqref{cpv} with a
time-dependent potential $V(t,x)$, that is
\begin{equation}\label{cpvpot}
\begin{cases}
\; \partial^2_t u + \Delta^2 u + V(t,x)u=F(t,x)\\
\;u(0,x)=u_0(x),\,\,
\partial_t u (0,x)=u_1(x).
\end{cases}
\end{equation}
Let us first introduce the notion of admissible pairs.
\begin{defi}
We say that the exponent pair $(q,r)$ is \emph{admissible} if
\[
2\leqslant q,r\leqslant \infty,\quad \frac{2}{q}+\frac{d}{r}
=\frac{d}{2},\quad d\geqslant 1 \quad (q,r,d)\neq (2,\infty,2).
\]
\end{defi}
We shall prove that the vibrating plate equation with a
time-dependent potential \eqref{cpvpot} is well-posed in
$\Hs^s(\R^d)\times\Hs^{s-2}(\R^d)$, with $s\in [2,d)$, for the
following class of potentials
\begin{tm}\label{teovib}
 Let $d\geq 3$, $I$ be either the interval $[0,T]$, $T>0$,  or $[0,+\infty)$, and assume $V(t,x)$ is a real valued potential such that
\begin{equation}\label{pot}
  V(t,x)\in L_I^{\alpha} \W^{s-2,\be} \quad\text{with}\quad \frac{2}{\alpha}+\frac{d}{\beta}=s+2,
\end{equation}
for some fixed $s\in[2,d),\,\alpha\in[1,\infty)$ and $\beta\in(1,\infty)$.
If $u_0\in \dot{H}^s,\,u_1\in \dot{H}^{s-2}$ and $F\in L_I^{\t q'}\W^{s-2,\t r'}$, for some admissible pair $({\t q},{\t r})$,
 then the Cauchy problem \eqref{cpvpot} has a unique solution $u\in C(I;\Hs^s(\R^d))\cap L^q(I;\W^{s,r}(\R^d))$, for all admissible pairs $(q,r)$, such that
 \[
 \| u\|_{L^q_I\W^{s,r}}\leq C_V \|u_0\|_{\Hs^{s}} +  C_V \|u_1\|_{\Hs^{s-2}} +  C_V \|F\|_{L^{\t q'}_I\W^{{s-2},{\t r'}}}.
\]
\end{tm}
The admissibility class of the potentials is represented in Figure \ref{fig2}. The corresponding problem for the
Schr\"odinger equation was studied in \cite{DaPiVi}, see Figure \ref{fig3}.

\begin{figure}[b]
\centering
\begin{minipage}{.45\textwidth}
\centering
 \begin{tikzpicture}[>=latex,scale=3]
  \draw[->] (-0.05,0) -- (1.2,0) node[below] {$1/\alpha$};
  \draw[->] (0,-0.1) -- (0,1.2) node[left] {$1/\beta$};
  \coordinate [label=left :\textcolor{black}{\small $1$}](a) at (0,1);
  \coordinate [label=below :\textcolor{black}{\small $1$}](b) at (1,0);
  \draw[very thin, dashed, color=gray]  (a) -- (1,1);
  \draw[very thin, dashed, color=gray]  (b) -- (1,1);

  \draw[thick, color=blue!45] (0,4/7) -- (1,2/7); 
  \coordinate[label=right:\textcolor{blue!45}{\small $s=2$}] (d) at (0.6,0.45);
  \fill[blue!45] (1,2/7) circle (0.5pt);
  \draw[blue!45] (0,4/7) circle (0.5pt);
  \draw[thick, color=blue!80!black] (0,6/7) -- (1,4/7); 
  \coordinate[label=right:\textcolor{blue!80!black}{\small $s=4$}] (e) at (0.65,0.75);
  \fill[blue!80!black] (1,4/7) circle (0.5pt);
  \draw[blue!80!black] (0,6/7) circle (0.5pt);

  \coordinate[label=above:\textcolor{black}{\small $V\in L_I^\alpha\W^{s-2,\beta}$}] (e1) at (0.6,1);
  \coordinate [label=left :\textcolor{black}{\small $2/d$}](c) at (0,2/7);
  \coordinate [label=left :\textcolor{black}{\small $4/d$}](c) at (0,4/7);
  \coordinate [label=left :\textcolor{black}{\small $6/d$}](c1) at (0,6/7);
  \draw[very thin, dashed, color=gray]  (0,2/7) -- (1,2/7);
  \draw[very thin, dashed, color=gray]  (0,4/7) -- (1,4/7);
  \draw[very thin, dashed, color=gray]  (0,6/7) -- (1,6/7);
\end{tikzpicture}
  \caption{Admissibility conditions on the potentials for the VP equation.}\label{fig2}
\end{minipage}
\hspace*{2mm}
\begin{minipage}{.45\textwidth}
\centering
\begin{tikzpicture}[>=latex,scale=3]
  \draw[->] (-0.05,0) -- (1.2,0) node[below] {$1/\alpha$};
  \draw[->] (0,-0.1) -- (0,1.2) node[left] {$1/\beta$};
  \coordinate [label=left :\textcolor{black}{\small $1$}](a) at (0,1);
  \coordinate [label=below :\textcolor{black}{\small $1$}](b) at (1,0);
  \draw[very thin, dashed, color=gray]  (a) -- (1,1);
  \draw[very thin, dashed, color=gray]  (b) -- (1,1);

  \draw[thick, color=red] (0,2/7) -- (b); 
  \fill[red] (b) circle (0.5pt);
  \draw[red] (0,2/7) circle (0.5pt);
  \coordinate[label=right:\textcolor{red}{\small Schr\"odinger}] (e) at (0.8,0.15);
  \draw[thick, color=blue] (0,4/7) -- (1,2/7); 
  \fill[blue] (1,2/7) circle (0.5pt);
  \draw[blue] (0,4/7) circle (0.5pt);
  \coordinate[label=right:\textcolor{blue}{\small vibrating plate}] (d) at (0.7,0.45);

  \coordinate[label=above:\textcolor{black}{\small $V\in L_I^\alpha L^\beta$}] (e2) at (0.6,1);
  \draw[very thin, dashed, color=gray]  (0,2/7) -- (1,2/7);
  \draw[very thin, dashed, color=gray]  (0,4/7) -- (1,4/7);
  \coordinate [label=left :\textcolor{black}{\small $2/d$}](c) at (0,2/7);
  \coordinate [label=left :\textcolor{black}{\small $4/d$}](c) at (0,4/7);
\end{tikzpicture}
\caption{Comparison between Schr\"odinger and VP equations ($s=2$).}\label{fig3}
\end{minipage}
\end{figure}
Differently from what has been achieved in the
literature so far (see, e.g.\cite{AgliardiCicognani,agza,corzuc}), where
the study of the Cauchy problem for the VP equation was only
\emph{local} in time, thanks to the Strichartz estimates, we are
able to obtain also a \emph{global} solution to the problem.

The proof of the theorem (given in Section \ref{secwell1})
requires an H\"older-type inequality for the homogeneous Sobolev
spaces. Indeed,  in Section \ref{seckp},  we shall use the
Kato-Ponce's inequality for homogeneous Sobolev spaces  and  the
Sobolev embeddings to obtain  the  H\"older-type inequality

\begin{equation*}
 \|f\cdot g\|_{\W^{s,r}}\lesssim
 \|f\|_{\W^{s_1,r_1}}\|g\|_{\W^{s_2,r_2}},
\end{equation*}
with $$\frac1r=\frac1{r_1}+\frac1{r_2}+\frac{s-s_1-s_2}{d},\quad
0\leq s\leq \min\{s_1,s_2\}, \quad r,r_1,r_2\in(1,\infty).$$

The  classical results on the existence of solutions to some
stationary  nonlinear equation (the so-called ground states)
\cite{beli} will play a central role in the construction of a
ground state $v$ for the stationary vibrating plate equation
\begin{equation*}
  \Delta^2 v - v + W(x)v = 0,
\end{equation*}
with  $W$ being a  potential in the Schwartz class $\cS(\R^d)$.
The ground state $v$ will be used to build up a suitable solution
to \eqref{cpvpot} which shall be employed to prove the sharpness of
\eqref{pot} (see Sections \ref{secwell2} and \ref{secwell3}).
\vskip0.1truecm
 We plan further investigations on this topic,
concerning the extension of  the results attained to homogeneous
Triebel-Lizorkin spaces and  the study of Strichartz estimates in
Wiener amalgam spaces for the LVP equation, by exploiting the
recent  results for the Schr\"odinger one in \cite{cornic}.

Further, observe that the Cauchy problem
\[
\begin{cases}
\; \partial^2_t u + \Delta^2 u + V(x)u=F(t,x)\\
\;u(0,x)=u_0(x),\,\,
\partial_t u (0,x)=u_1(x).
\end{cases}
\]
with the time independent potentials $V(x)$, is not allowed in
Theorem \ref{teovib} (since $\alpha<\infty$).
 However, this topic for the Schr\"odinger propagator has been
extensively studied in the literature. For instance, by combining
the arguments in \cite{rodshc} with the counterexamples given in
\cite{govevi}, it is possible to prove that the decay
$|x|^{-2-\epsilon}$ for $\epsilon>0$ and $|x|>>1$ is the sharp
decay to be required to the potential $V(x)$, in order to
guarantee the validity of Strichartz estimates for the
corresponding perturbed Schr\"odinger propagator
$e^{it\left(\Delta+V(x)\right)}$.

In view of Theorem \ref{teovib}, taking formally $\alpha= \infty$,
it seems that the natural (and eventually sharp) decay to be
imposed on $V(x)$ is $|x|^{-4-\epsilon}$ in order to guarantee
Strichartz estimates (for instance similar to the ones in Theorem
\ref{teovib} with s = 2). This topic will be also developed in the
future.

\vskip0.5truecm
 We conclude  by fixing some notations. For $1\leq
p\leq\infty$, let $p'$ be the conjugate exponent of $p$
($1/p+1/p'=1$). We shall use  $A\lesssim B$ to mean that there
exists a constant $c>0$ such that $A\leq cB$ and $A\asymp B$ means
that $A\lesssim B \lesssim A$.   For any subinterval $I$ of $\R$
(bounded or unbounded), we define the mixed space-time norms
\[
 \|u\|_{L_I^q X} := \Big(\int_I \|u(t,\cdot)\|_X dt\Big)^{1/q}
\]
with $X$ being either $L_{\phantom{I}}^r(\Rd)$ or $\W_{\phantom
I}^{s,r}(\Rd)$, and, similarly, for the spaces $C_I X$. When
$I=[0,+\infty)$ we write simply $L^qX$ in place of $L^q_IX$ and,
similarly, when $L^q_I$ is replaced by $C_I$. The space $\cS(\Rd)$
denotes the Schwartz class and $\cS'(\R^d)$ is its dual (the space
of tempered distributions). We write  $xy=x\cdot y$, $x,y\in
\R^d$, for the scalar product on $\R^d$. The Fourier transform is
given by ${\wh{f}}(\xi)=\F f(\xi)=\int e^{-i t\xi}f(t)dt$.

\section{Kato-Ponce's inequality}\label{seckp}

In this section we shall present the   Kato-Ponce's inequality for
the homogeneous Sobolev spaces. We first recall  the definitions
of (potential) Sobolev spaces.

\subsection{Sobolev spaces}
Let $s\in \R,\,1\leqslant r\leqslant\infty$ and $f\in\cS'(\Rd)$.
The (potential) \emph{Sobolev spaces} $W^{s,r}$ and the
(potential) \emph{homogeneous Sobolev spaces} $\W^{s,r}$ are
defined by
\[
\begin{split}
 W^{s,r}(\Rd)=&\Big\{f\in\cS'(\Rd): \|f\|_{W^{s,r}}=\|\l\cD\r^s f \|_{L^r}<\infty\Big\}\\
 \W^{s,r}(\Rd)=&\Big\{f\in\cS'(\Rd): \|f\|_{\W^{s,r}}:=\||\cD|^s f
 \|_{L^r}<\infty\Big\},
\end{split}
\]
where the fractional differentiation operators
$\l\cD\r^s$ and $|\cD|^s$ are the Fourier multipliers defined by
\begin{equation}\label{mult}
 \wh{\l\cD\r^s f}(\xi):= \l \xi\r^s \wh f(\xi) \quad \text{ and }\quad \wh{|\cD|^s f}(\xi):= |\xi|^s \wh f(\xi).
\end{equation}
In particular, if $s=2$ then $\l\cD\r^2= I-\De$,  where $I$ is the
identity operator, and $|\cD|^2=-\De$. The operators
$\l\cD\r^{-s}$ with $s>0$ and $|\cD|^{-s}$ with $0<s<d$ are also
known as the \emph{Bessel} and \emph{Riesz potentials} of order s,
respectively (see, e.g., \cite{stein}). If $r=2$ these spaces are
also denoted by $H^s(\Rd)$ and $\Hs^s(\Rd)$. For $s\geq 0$,
$W^{s,r}$ and its homogeneous counterpart $\W^{s,r}$ are related
as follows
\[
W^{s,r}=L^r\cap \W^{s,r} ,\quad \text{with} \quad
\|\cdot\|_{W^{s,r}}\asymp\|\cdot\|_{L^r}+\|\cdot\|_{\W^{s,r}}.
\]
Recall that the homogeneous Sobolev space $\W^{s,r}$ is a
seminormed space and $\|f\|_{\W^{s,r}}=0$ if ad only if $f$ is a
polynomial.

Consider  the dilation operator $S_\la (f)(x)= f(\la x)$ for
$\la>0$, then the dilation properties for homogeneous Sobolev
spaces read
\begin{equation*}
\|S_\lambda f\|_{\W^{s,r}}=\lambda^{s-d/r}\|f\|_{\W^{s,r}},\quad
s\in \R, \,\,1\leq r\leq \infty.
\end{equation*}
Finally, the following \emph{homogeneous Sobolev embeddings}
\cite{BerLof}  will be useful  in the sequel.
\begin{lemma}\label{teoemb}
Assume $f\in \cS(\R^d),\, s\in\R$ and $1<r<\infty$, then
\[
 \|f\|_{\W^{s,r}}\lesssim \|f\|_{\W^{s_1,r_1}}
\]
where $s\leq s_1$ and $1<r_1\leqslant r<\infty$ are such that
$s-d/r=s_1-d/r_1$.
\end{lemma}

For further  properties of these spaces we address the reader to,
e.g.,  \cite{BerLof}.

\subsection{Kato-Ponce's inequality}
The \emph{Kato-Ponce's inequality}, established in \cite{ponce} by
combining the original argument in \cite{katoponce} with the
general version of the Coifman and Meyer's result in
\cite{CoifMey}, reads as follows
\begin{tm}\label{teokp}
 Suppose $f,g\in \cS(\R^d),\, s\geq 0$ and $1<r<\infty$, then
\begin{equation}\label{eqkp}
 \|\l\cD\r^s(f\cdot g)\|_{L^r}\lesssim \|f\|_{L^{r_1}} \|\l\cD\r^s g\|_{L^{r_2}} + \|\l\cD\r^s f\|_{L^{r_3}} \| g\|_{L^{r_4}}
\end{equation}
with
\begin{equation}\label{indrel}\frac1r=\frac1{r_1}+\frac1{r_2}=\frac1{r_3} +
\frac1{r_4},\quad r_2,\,r_3\in (1,\infty).\end{equation}
\end{tm}
Observe that this implies $r_1,r_4\in (1,\infty]$. Combining
dilation properties of the operator $S_\la$ with the Lebesgue's
dominated convergence Theorem,  it follows
corresponding inequality for the homogeneous Sobolev spaces.
\begin{tm}\label{teokphom}
 Suppose $f,g\in \cS(\R^d),\, s\geq 0$ and $1<r<\infty$, then
\begin{equation}\label{eqkph}
 \||\cD|^s(f\cdot g)\|_{L^r}\lesssim \|f\|_{L^{r_1}} \||\cD|^s g\|_{L^{r_2}} + \||\cD|^s f\|_{L^{r_3}} \| g\|_{L^{r_4}}
\end{equation}
with the indices' relations \eqref{indrel}.
\begin{proof}
The inequality is well-known. However, we detail the proof for the sake of clarity. We apply the Kato-Ponce's inequality \eqref{eqkp} to the rescaled
product $S_\la(f\cdot g)=S_\la(f)S_\la( g)$
\begin{equation}\label{eqdim}
   \|\l\cD\r^s( S_\la(f\cdot g))\|_{L^r}\lesssim \| S_\la f\|_{L^{r_1}} \|\l\cD\r^s S_\la g\|_{L^{r_2}} + \|\l\cD\r^s  S_\la f\|_{L^{r_3}}
   \|  S_\la g\|_{L^{r_4}},
\end{equation}
where $r_2,\,r_3\in (1,\infty)$ are such that
$1/r=1/r_1+1/r_2=1/r_3 + 1/r_4$. Combining the commutativity
property
\[
 \l\cD\r^s S_\la (f)(x) = \la^s S_\la\big((\la^{-2}-|D|^2)^{s/2}f\big)(x),\qquad \la>0
\]
with $\|S_\la f\|_{L^r}=\la^{-d/r}\|f\|_{L^r},\; \la>0$, the
equation \eqref{eqdim} can be rewritten as
\[
 \|(\la^{-2}-|D|^2)^{s/2}(f\cdot g)\|_{L^r} \lesssim \|f\|_{L^{r_1}} \|(\la^{-2}-|D|^2)^{s/2} g\|_{L^{r_2}} + \|(\la^{-2}-|D|^2)^{s/2}
 f\|_{L^{r_3}} \| g\|_{L^{r_4}}.
\]
We reach the claim if we show that
\[
 \lim_{\la\to+\infty} \|(\la^{-2}-|D|^2)^{s/2}f\|_{L^r} = \||D|^s f\|_{L^r}.
\]
But this fact is simply a consequence of the Lebesgue's dominated
convergence Theorem applied to the functions
\[
 f_\la:= |(\la^{-2}-|D|^2)^{s/2}f|^r. \qedhere
\]
\end{proof}
\end{tm}

\begin{rem}
In the case $d=2$ the above inequality \eqref{eqkph} was already
mentioned in \cite{MuPiTaoTh}.
\end{rem}

From this result we derive the following H\"older-type inequality
\begin{cor}
 Suppose $f,g\in \cS(\R^d),\, s\geq 0$ and $1<r<\infty$, then
\begin{equation}\label{kpholder}
 \|f\cdot g\|_{\W^{s,r}}\lesssim \|f\|_{\W^{s_1,r_1}}\|g\|_{\W^{s_2,r_2}}
\end{equation}
with
\begin{equation}\label{indrel2}\frac1r=\frac1{r_1}+\frac1{r_2}+\frac{s-s_1-s_2}{d},\quad
0\leq s\leq \min\{s_1,s_2\}, \quad
r,r_1,r_2\in(1,\infty).\end{equation}
\begin{proof}
Using the estimate  \eqref{eqkph} and the embeddings of Lemma
\ref{teoemb} we attain the desired estimate.
\end{proof}
\end{cor}

\begin{rem}
 The Kato-Ponce's inequality in the case $s=s_1=s_2=0$ reduces to the classical H\"older's inequality.
\end{rem}

\section{Strichartz estimates}\label{ST}
In this section, we first recall the well-known Strichartz
estimates for the Schr\"odinger equation. Then, we use them to
derive the corresponding ones for the Cauchy problem \eqref{cpv}.

\subsection{Schr\"odinger equation}\label{secsch}
Consider the linear  Schr\"odinger equation
\begin{equation}\label{scho}
 i\partial_t u(t,x)+\De u(t,x)=F(t,x),
\end{equation}
with initial data $u(0,x)=u_0(x)$. The solution $u(t,x)$ can be
formally written in the integral form
\[
 u(t,\cdot)=e^{it\De}u_0+\int_0^t e^{i(t-s)\Delta} F(s)\,ds,
\]
where, for every fixed $t$,  the Schr\"odinger
propagator $e^{it\De}$ is a Fourier multiplier  with symbol $e^{-i t |\xi|^2}$, $\xi\in\R^d$.
We recall that, given a function $\sigma$ on $\R^d$ (the so-called
symbol of the multiplier), the corresponding Fourier multiplier
operator $H_\sigma$ is formally defined by
$$H_\sigma f (x)= (2\pi )^{-d}\int_{\R^d} e^{i x \xi} \sigma(\xi)
\hat{f} (\xi) \,d\xi.
$$

So, continuity properties for multipliers in suitable spaces yield
estimates for the solution to the equation \eqref{scho}.
From the previous formula, when $F(t,x)=0$,  the following
\emph{fixed-time estimates} can be obtained
\begin{equation}\label{displp}
\|e^{it\Delta}u_0 \|_{L^{r}}\lesssim
|t|^{-d\left(\frac{1}{2}-\frac{1}{r}\right)}\|u_0\|_{L^{r'}},\quad
2\leqslant r\leqslant\infty.
\end{equation}
Indeed, if the initial data $u_0$ displays a suitable
integrability in space, then the evolution will have a power-type
decay in time. In the particular case $r=2$, we have  the
\emph{conservation law of  energy}
$\|e^{it\Delta}u_0\|_{L^{2}}=\|u_0\|_{L^{2}}.$

By combining the above dispersive estimates with some duality
arguments, one can obtain the celebrated Strichartz estimates
\cite{cazenave,gvelo,ktao}.
\begin{tm}
 Let $I\subseteq \R$, $(q,r)$ and $(\t q,\t r)$ be admissible pairs. Then, for all initial data $u_0\in L^2(\R^d)$ and forcing term
  $F\in L_I^{\t q'}L_{\phantom{I}}^{ \t r'}$, we have the homogeneous Strichartz estimates
 \begin{equation}\label{strihom}
 \|e^{it\De}u_0\|_{L_I^{q}L^{r}} \lesssim \|u_0\|_{L^2}
\end{equation}
and the inhomogenous Strichartz estimates
\begin{equation}\label{strinhom}
 \Big\|\int_0^t e^{i(t-s)\Delta} F(s) ds\Big\|_{L_I^qL_{\phantom{I}}^{ r}} \lesssim \|F\|_{L_I^{\t q'}L_{\phantom{I}}^{ \t
 r'}},
\end{equation}
with constants that do not depend on the length of $I$.
\end{tm}

\subsection{Vibrating plate equation}

In order to find a solution to the homogeneous VP equation \eqref{cpv}, (e.g., with $F=0$),  we compute the
Fourier transform of the problem \eqref{cpvom}, obtaining the
following differential equation
\begin{equation}\label{cpvtrasf}
\begin{cases}
\; \partial^2_t \wh u + |\xi|^4 \wh u=0\\
\;\wh u(0,\xi)=\wh u_0(\xi),\,\,
\partial_t \wh u (0,\xi)=\wh u_1(\xi).
\end{cases}
\end{equation}
It is straightforward to solve it and, taking the inverse Fourier
transform of the solution, we can express the solution to
\eqref{cpvom} in the form
\begin{equation}\label{solhom}
 u(t,\cdot) = K'(t) u_0 + K(t) u_1
\end{equation}
where,
\begin{equation}\label{solhom2}
 K'(t) = \cos(t\De), \quad\quad K(t)= \frac{\sin(t\De)}{\De}.
\end{equation}
 Here, for every fixed $t$, the propagators $K(t), K'(t)$ are Fourier multipliers  with symbols $\cos(t|\xi|^2)$,
 $\sin(t|\xi|^2)/|\xi|^2,\,\xi\in\Rd$.

Next, by Duhamel's formula, the solution to the inhomogeneous
equation \eqref{cpv} can be formally written in the integral form
\begin{equation}
u(t,\cdot)=\cos(t\De)u_0 + \frac{\sin(t\De)}{\De}u_1 + \int_0^t
\frac{\sin((t-s)\De)}{\De} F(s)\, ds.
\end{equation}
The simply but efficacious formulas
\[
 \cos(t\De)=\frac{e^{it\De}+e^{-it\De}}{2}, \qquad \frac{\sin(t\De)}{\De} = \frac{e^{it\De}-e^{-it\De}}{2i \De}
\]
show that estimates for $K'$ follow directly from the ones for the
Schr\"odinger equation, whereas estimates on $K$ can be obtained
by the continuity properties of the propagator
$\displaystyle{\frac{e^{it\De}}{\De}}$.

 The commutativity
property of the Fourier multipliers $e^{it\De}$ and $\De$,
combined with \eqref{displp} yields the \emph{fixed-time
estimates}
\[
  \|e^{it\De} u_0\|_{\W^{s,r}}\lesssim
  |t|^{-d(\frac{1}{2}-\frac{1}{r})}\|u_0\|_{\W^{s,r'}},\,\quad\,
 \Big\|\frac{e^{it\De}}{\De} u_1\Big\|_{\W^{s,r}}\lesssim
 |t|^{-d(\frac{1}{2}-\frac{1}{r})}\|u_1\|_{\W^{s-2,r'}},
\]
with $2\leqslant r\leqslant\infty,\, s\in \R$. Moreover, the
sharpness of  these estimates follows from the optimality of
\eqref{displp}. Then, using the  Strichartz estimates for the
Schr\"odinger equation and the commutativity of Fourier
multipliers we obtain the following  estimates
\begin{tm}\label{viteostri}
Let $d\geqslant 1$ and $s\in\R$. Then, the following homogeneous
Strichartz estimates
\begin{equation}\label{vistrihom}
   \|e^{it\De} u_0\|_{L^q_I\W^{s,r}}\lesssim
  \|u_0\|_{\Hs^{s}},\,\quad \Big\|\frac{e^{it\De}}{\De}u_1\Big\|_{L_I^q\W_{\phantom{I}}^{s,r}}\lesssim \|u_1\|_{\Hs^{s-2}}
\end{equation}
and inhomogeneous Strichartz estimates
\begin{equation}\label{vistrinhom}
\Big|\Big|\int_{0}^t \frac{e^{i(t-s)\Delta}}{\De}F(s)\,ds\Big|\Big|_{L_I^q \W_{\phantom{I}}^{s,r}}\lesssim ||F||_{L_I^{\tilde{q}'}\W_{\phantom{I}}^{s-2,\tilde{r}'}}.
\end{equation}
hold for any admissible pairs $(q,r)$ and $(\t q,\t r)$.
\begin{proof}
Estimates \eqref{vistrihom} and  \eqref{vistrinhom} follow
 from the the commutativity property of the Fourier
multipliers $e^{it\De}$ and $\De$ and the Strichartz estimates
\eqref{strihom} and \eqref{strinhom}, respectively.
\end{proof}
\end{tm}
Observe that the admissible pairs for the VP equation are the
admissible Schr\"odinger couples.
\begin{cor}[Strichartz estimates for the vibrating plate]\label{corstri}
 Let $d\geq 1$, $s\in \R$ and $(q,r),\,(\t q,\t r)$ be admissible pairs. If $u$ is a  solution to the Cauchy problem \eqref{cpv}, then
\begin{equation}\label{vistrisolhom}
 \|u\|_{L_I^q \W_{\phantom{I}}^{s,r}}\lesssim \|u_0\|_{\Hs^{s}} + \|u_1\|_{\Hs^{s-2}} + \|F\|_{L_I^{\t q'}\W_{\phantom{I}}^{s-2,\t r'}}.
\end{equation}
\end{cor}

\section{Application to vibrating plate equations with time-dependent potentials}\label{secwell}
The Strichartz estimates obtained above  are a fundamental tool to
obtain the well-posedness of the Cauchy problem \eqref{cpvpot}.
First, we present the proof of the main theorem stated in the
introduction. Secondly, we prove the sharpness of the
admissibility conditions for the potentials in \eqref{pot}.  The
basic idea  is to use  suitable rescaling arguments for a standing
wave, solution to  the homogeneous equation
\begin{equation}\label{eqomo}
 \partial^2_t u + \Delta^2 u + V(t,x)u=0.\\
\end{equation}

\subsection{Proof of the Theorem \ref{teovib}}\label{secwell1}
\begin{proof}
The proof uses the arguments of \cite[Theorem 1.1]{DaPiVi} (see also \cite[Theorem 6.1]{cornic}), based on the Banach-Caccioppoli contraction theorem.  Strichartz estimates represent the key to achieve the contraction.

First,  we need to  choose the suitable space and the mapping  to apply the contraction theorem.
Let $J=[0,\delta]$ be a small interval and, for any admissible pair $(q,r)$, $s\in [2,d)$, set
\[
 Z_{q,r}^s=L^q(J;\W^{s,r}(\R^d))
\]
and
\[
 Z=C(J;\Hs_{\phantom J}^s(\R^d))\cap Z_{2,2d/(d-2)}^s \quad\text{with}\quad \|v\|_Z=\max\bigg\{\|v\|_{L_J^{\infty}\Hs_{\phantom{J}}^s}\;,\;\|v\|_{Z_{2,2d/(d-2)}^s}\bigg\}.
\]
Notice that, by complex interpolation (see \cite[Parag. 5.2.5]{trieb} or  \cite[Parag. 6.4]{BerLof}),
\[
\big[C(J;\Hs_{\phantom J}^s(\R^d)) \;;\; Z_{2,2d/(d-2)}^s\big]_{[\th]}=Z_{q,r}^s \quad \text{with }\frac{1}{q}=\frac{\th}{2} \text{ and }\frac{1}{r}=\frac{1-\th}{2} +\frac{\th(d-2)}{2d}.
\]
This gives $2/q+d/r=d/2$, i.e., the admissibility condition;  so the space Z is embedded in all admissible spaces $Z_{q,r}^s$.

Consider now, for any $v\in Z$, the mapping
\[
 \phi(v)=\cos(t\De)u_0+\frac{\sin(t\De)}{\De} u_1 +\int_0^t \frac{\sin(t\De)}{\De} \big[F(s) - V(s)v(s)\big]ds.
\]
From Corollary \ref{corstri} it is easy to see that the subsequent estimates hold:
\begin{equation}\label{eqdim2}
\|\phi(v)\|_{Z_{q,r}^s}\leqslant C_0 \|u_0\|_{\Hs^{s}} + C_0\|u_1\|_{\Hs^{s-2}} + C_0\|F\|_{Z_{\t q',\t r'}^{s-2}} + C_0\|Vv\|_{Z_{q'_0,r'_0}^{s-2}}
\end{equation}
for all $(q,r),(\t q,\t r)$ and $(q_0, r_0)$ admissible.
In order to estimate the product $Vv$ appropriately we shall prove the case $2\leqslant\alpha<\infty$ and $1\leqslant\alpha< 2$ separately.

Starting with $2\leqslant\alpha<\infty$, among such admissible pairs, there is a pair $(q_0,r_0)$ such that
\[
 \frac{1}{q_0} = \frac{1}{2}-\frac{1}{\alpha} \quad \text{and}\quad \frac{1}{r_0}=-\frac{1}{\be}+\frac{1}{2}+\frac{s+1}{d}.
\]
In fact, from  $\alpha\geqslant 2$ and \eqref{pot}, it follows $q_0\geqslant 2$ and
\[
\frac{2}{q_0}+\frac{d}{r_0}=2\Big(\frac{1}{2}-\frac{1}{\alpha}\Big)+d\Big(-\frac{1}{\be}+\frac{1}{2}+\frac{s+1}{d}\Big)=\frac{d}{2}.
\]
We choose such a pair in \eqref{eqdim2}. Then, we apply Holder's inequality and the estimate \eqref{kpholder}, obtaining
\[
 \|\phi(v)\|_{Z_{q,r}^s}\leqslant C_0\|u_0\|_{\Hs^{s}} + C_0\|u_1\|_{\Hs^{s-2}} + C_0\|F\|_{Z_{\t q',\t r'}^{s-2}} + C_0\|V\|_{Z_{\alpha,\be}^{s-2}}\|v\|_{Z_{2,2d/(d-2)}^s}.
\]
By taking $(q,r)=(2,2d/(d-2))$ and $(q,r)=(\infty,2)$ one deduces that $\phi: Z\to Z$ (the fact that $\phi(u)$ is \emph{continuous} in $t$ when valued in $\Hs^s$ follows from a classical limiting
argument, see \cite[Theorem 1.1, Remark 1.3]{DaPiVi}). Also, since $\alpha<\infty$, if $J$ is small enough, we get $C_0\|V\|_{L_J^{\al}\W^{s-2,\be}}<1/2$, hence  $\phi$ is a contraction.
The contraction theorem then states the existence and the uniqueness of the  solution in $J$. By iterating this argument a finite number of times one obtains a solution in $[0,T]$ if $T<\infty$ or in $\R$ if $T=\infty$ (again, see \cite[Theorem 1.1]{DaPiVi}).

In the case $1\leqslant\alpha<2$, there is a pair $(q_0,r_0)$ such that
\[
 \frac{1}{q_0} = 1-  \frac{1}{\alpha} \quad \text{and}\quad \frac{1}{r_0}=-\frac{1}{\be}+\frac{1}{2}+\frac{s}{d}.
\]
In fact, from $1\leqslant\al<2$ and \eqref{pot}, it follows $q_0> 2$ and

\[
\frac{2}{q_0}+\frac{d}{r_0}=2\Big(1-\frac{1}{\alpha}\Big)+d\Big(-\frac{1}{\be}+\frac{1}{2}+\frac{s}{d}\Big)=\frac{d}{2}.
\]
We choose such a pair in \eqref{eqdim2}, in order to apply H\"older's inequality and the estimate \eqref{kpholder} again. We obtain
\[
 \|\phi(v)\|_{Z_{q,r}^s}\leqslant C_0 \|u_0\|_{\Hs^{s}} + C_0\|u_1\|_{\Hs^{s-2}} + C_0\|F\|_{Z_{\t q',\t r'}^{s-2}} +  C_0\|V\|_{Z_{\alpha,\be}^{s-2}}\|v\|_{Z_{\infty,2}^s}.
\]
Finally, the same arguments as for the previous case yield  the desired result.
\end{proof}

\begin{rem}
Let us make a few comments about the  constraints on the indices $d,s,\alpha, \beta$. They are mainly due to the need
of applying the Kato-Ponce's estimate to the product $\| V v\|$ in \eqref{eqdim2}.
First, we must exclude the dimension $d=2$, since it forces $s=2, \alpha=\beta=1$ and Kato-Ponce's inequality does not hold.
 So, we start with the dimension $d\geq 3$. Similarly, Kato-Ponce does not hold in the cases: (i) $s<2$, (ii) $\beta=1$.
If $s\geqslant d$ there are no pairs ($\alpha,\beta)\,\in [1,\infty)\times(1,\infty)$ that satisfy \eqref{pot} and the same holds when $\beta=\infty$, for all $\al$.
Finally, observe that the condition on $\alpha\in [1,\infty)$, forces $\beta$ to live in the bounded range $\big(d/(s+2),d/s\big]$, see Figure \ref{fig1} and also Figure \ref{fig2}.
\end{rem}
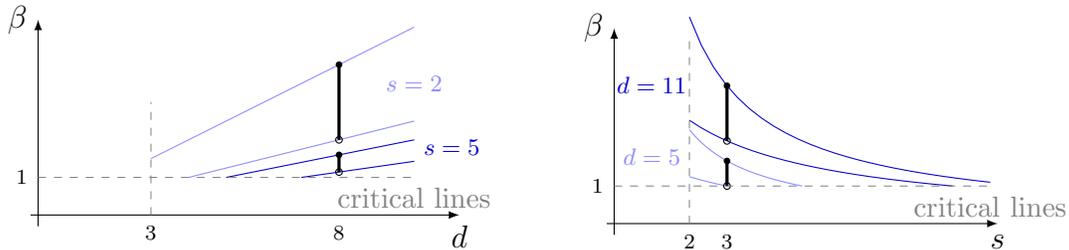
\begin{figure}[b]
\centering
\begin{tikzpicture}[>=latex, domain=3:10, scale=0.5]
  \draw[->] (-0.2,0) -- (11.2,0) node[below] {$d$};
  \draw[->] (0,-0.3) -- (0,5.2) node[left] {$\beta$};
  \coordinate [label=below :\textcolor{black}{$\scriptstyle 3$}](a) at (3,0);
  \coordinate [label=below :\textcolor{black}{$\scriptstyle 8$}](b) at (8,0);
  \coordinate [label=left :\textcolor{black}{$\scriptstyle 1$}](c) at (0,1);

  \draw[color=blue!45, domain=4:10]   plot (\x,{\x/(4)}); 
  \draw[color=blue!45]   plot (\x,{(\x)/(2)})   ; 
  \coordinate [label=below :\textcolor{blue!45}{\scriptsize $s=2$}](d) at (10,4);
  \draw[color=blue!80!black, domain=7:10]   plot (\x,{\x/(7)}); 
  \draw[color=blue!80!black, domain=5:10]   plot (\x,{(\x)/(5)}); 
  \coordinate [label=below :\textcolor{blue!80!black}{\scriptsize $s=5$}](e) at (11,2.3);

  \draw[very thick, color=black]  (8,4) -- (8,2)    ;
  \fill[black] (8,4) circle (2.5pt);
  \draw[black] (8,2) circle (2.5pt);
  \draw[very thick, color=black]  (8,8/7) -- (8,8/5)    ;
  \fill[black] (8,8/5) circle (2.5pt);
  \draw[black] (8,8/7) circle (2.5pt);

  \draw[very thin, dashed, color=gray, domain=0:10] plot (\x,1) node[below] {\small critical lines};
  \draw[very thin, dashed, color=gray] (3,0)--(3,3);
\end{tikzpicture}
\qquad
\begin{tikzpicture}[>=latex, domain=2:10, scale=0.5]
  \draw[->] (-0.2,0) -- (10.2,0) node[below] {$s$};
  \draw[->] (0,-0.3) -- (0,5.2) node[left] {$\beta$};
  \coordinate [label=below :\textcolor{black}{$\scriptstyle 2$}](f) at (2,0);
  \coordinate [label=below :\textcolor{black}{$\scriptstyle 3$}](g) at (3,0);
  \coordinate [label=left :\textcolor{black}{$\scriptstyle 1$}](h) at (0,1);

  \draw[color=blue!45!,domain=2:5]     plot (\x,{(5)/(\x)}) ; 
  \draw[color=blue!45!,domain=2:3]     plot (\x,{(5)/(\x+2)}); 
  \coordinate [label=below :\textcolor{blue!45!}{\scriptsize $d=5$}](i) at (1,2.3);
  \draw[color=blue!80!black!,domain=2:10]    plot (\x,{(11)/(\x)}); 
  \draw[color=blue!80!black,domain=2:9]    plot (\x,{(11)/(\x+2)}); 
  \coordinate [label=below :\textcolor{blue!80!black}{\scriptsize $d=11$}](l) at (1,4.2);

  \draw[very thick, color=black]  (3,5/3) -- (3,1);
  \fill[black] (3,5/3) circle (2.5pt);
  \draw[black] (3,1) circle (2.5pt);
  \draw[very thick, color=black]  (3,11/3) -- (3,11/5);
  \fill[black] (3,11/3) circle (2.5pt);
  \draw[black] (3,11/5) circle (2.5pt);

  \draw[very thin, dashed, color=gray, domain=0:10] plot (\x,1) node[below] {\small critical lines};
  \draw[very thin, dashed, color=gray] (2,0)--(2,5);
\end{tikzpicture}
\caption{The black vertical segment shows the admissibility range for the exponent $\beta$. For every fixed dimension $d\geq 3$,
if the regularity $s$ increases the range decreases. On the other hand,  for every fixed regularity $s$,
if the dimension $d$ grows the range increases as well.}\label{fig1}
\end{figure}

In the limit case $s=2$ condition \eqref{pot} becomes
\begin{equation}\label{eqs2}
 V(t,x)\in L_I^\al L^\be \quad\text{with}\quad \frac{2}{\al}+\frac{d}{\be}=4,
\end{equation}
for some fixed $\alpha\in [1,\infty)$ and $\beta\in (d/4,d/2]$. In this case the potentials are reduced to the usual Lebesgue
spaces in the space variable, and we can compare the results for the VP equation with the corresponding ones  for Schr\"odinger
equation, exhibited in  \cite[Theorem 1.1]{DaPiVi}. Indeed, the class of admissible potentials for the Schr\"odinger equation is given by
\[
  V(t,x)\in L_I^\al L^\be \quad\text{with}\quad \frac{2}{\al}+\frac{d}{\be}=2,
\]
for some fixed $\alpha\in [1,\infty)$ and $\beta\in (d/2,\infty]$.
This means that the effect of composing two Schr\"odinger-type
propagators reduces notably the class of admissible potentials:
the exponent  $\beta$ goes from the unbounded region
$(d/2,\infty]$ to the bounded one $(d/4,d/2]$. Considering the
reciprocal $1/\beta$, we see that the class of admissible
potentials is shifted by a factor $2/d$, as shown in Figure
\ref{fig3}.

\begin{rem}
If $I=[0,T]$ is a bounded set, by H\"older's inequality in time, assumption \eqref{pot} can be relaxed to
\[
 V(t,x)\in L_I^{\alpha} \W^{s-2,\be} \quad\text{with}\quad \frac{2}{\alpha}+\frac{d}{\beta}\leqslant s+2.
\]
In this case the exponent $\beta$ is free to stay on the entire unbounded range $\big(d/(s+2),\infty\big)$.
\end{rem}

\begin{rem}\label{remresc}
Condition \eqref{pot} is quite natural in view of the following argument: the standard rescaling
$u_{\epsilon}(t,x)=u(\epsilon^2 t, \epsilon x)$ takes equation \eqref{cpvpot} into the equation
\[
 \partial_t u_\epsilon + \De^2 u_\epsilon + V_\epsilon(t,x) u_\epsilon=F(\epsilon^2 t,\epsilon x), \qquad V_\epsilon(t,x)=\epsilon^4 V(\epsilon^2 t, \epsilon x)
\]
and we have
\[
 \|V_\epsilon\|_{L^\alpha \W^{s-2,\beta}}=\epsilon^{-\frac{2}{\alpha}-\frac{d}{\be}+s+2}\|V\|_{L^\alpha \W^{s-2,\beta}}
\]
so that the $L^\alpha \W^{s-2,\beta}$ norm of $V_\epsilon$ is
independent of $\epsilon$ precisely when $\alpha, \be$ satisfy
\eqref{pot}. In the following two subsections we shall prove that
condition \eqref{pot} is compulsory for the global Strichartz
estimate to be true.  The most delicate part will be the
construction of a standing wave for the equation \eqref{eqomo}.
Hence, we now focus on this topic.
\end{rem}

\subsection{Existence of a ground state}\label{secwell2}

We shall prove the existence of a solution for the corresponding stationary equation of \eqref{eqomo}, that is
\begin{equation}\label{eqstat}
  \Delta^2 v - v + W(x)v = 0,
\end{equation}
for some potential $W(x)$. Even though the one dimensional case is
not covered by our results, it is useful to start from it, in
order to understand and solve the problem for any dimension. So,
for the moment, let $d=1$ and consider the \emph{Newton equation}
\begin{equation}\label{eqnew}
 -v''+ v - v^2 = 0.
\end{equation}
This classical equation comes out in the examination of solitary
waves for the Klein-Gordon and the Schr\"odinger equations. In
this specific example the explicit form of the solution is
well-known, precisely
\begin{equation}\label{sol1}
v(x)=\frac{3/2}{\cosh^2 (x/2)},  \qquad \text{where} \quad \cosh t
= \frac{e^t + e^{-t}}{2}.
\end{equation}
Notice that $v(x)>0$, for every $x\in\R$. Now, differentiating two
times this equation, one obtains
\[
 -v^{iv}+v''-2(v')^2-2vv''=0.
\]
The key idea is  to replace the second derivative $v''$ with its
expression in \eqref{eqnew}, so that
\[
v^{iv}-v+3v^2-2v^3+2(v')^2=0.
\]
Finally, using the explicit form of the solution $v$ in
\eqref{sol1}, we can compute  $(v')^2/v$ and show  the existence
of a solution $v$  for the following fourth order equation
\begin{equation}\label{eqfour}
 v^{iv}-v+W(x)v=0 \qquad  \text{where} \quad W(x)=-\frac{10}{3}v^2+5v.
\end{equation}
In this particular example the solution is known from the
beginning, given by  \eqref{sol1}, and we determine the potential
$W(x)$ accordingly. Let us underline  some  properties of this
solution we would like to have in any dimension $d$: \emph{i)}
$v>0$; \emph{ii)} $v$ is even: $v(x)=v(-x)$; \emph{iii)} $v\in
C^\infty(\R)$; \emph{iv)} $v$ and all its derivatives have
exponential decay at $\pm\infty$ (in particular,  $v\in\cS(\R)$).
Moreover,  we observe that $W(x)$ belongs to the Schwartz class
$\cS(\R)$. \vskip0.1truecm

Let us come  to the $d$-dimensional case ($d\geqslant 3$).
Inspired by this simple  example we shall prove the existence of a
solution to \eqref{eqstat}, arguing as follows: \emph{i)} prove
the existence of a solution for a well-known second order PDE,
\emph{ii)} apply the operator $\Delta$ to it and use  the original
second order PDE to eliminate remaining second order terms,
\emph{iii)} study the regularity and the decay of the solution.
The natural generalization from dimension $1$ to $d$ would be to
substitute the  power nonlinearity $v^2$ by $|v|^{p-1} v$, with
$p>1$. But this does not work for any dimension, as we shall see
presently. We shall make use of the following  classical result by
Berestycki and Lions \cite{beli}.
\begin{tm}\label{teobl}
 Suppose $d\geqslant 3$ and $g:\R\to\R$ be an odd continuous function such that $g(0)=0$. If the function $g$ satisfies the following conditions:
\begin{equation}\label{eqbl1}
 -\infty<\liminf_{s\to 0^+} \frac{g(s)}{s}\leqslant \limsup_{s\to0^+}\frac {g(s)}{s}=-m<0,
\end{equation}
\begin{equation}\label{eqbl2}
 -\infty\leqslant \limsup_{s\to+\infty}\frac {g(s)}{s^l}\leqslant 0,\quad \text{where}\quad l=\frac{d+2}{d-2},
\end{equation}
\begin{equation}\label{eqbl3}
\text{there exists } \zeta>0 \text{ such that } G(\zeta)=\int_0^\zeta g(s)ds>0,
\end{equation}
then the problem
\begin{equation}\label{prob}
 -\Delta v = g(v) \text{ in } \R^d,\qquad u\in H^1(\R^d),
\end{equation}
possesses a non trivial solution $v$ such that
\begin{itemize}
 \item[i)] $v>0$ on $\R^d$.
 \item[ii)] $v$ is spherically symmetric: $v(x)=v(r)$, where $r=|x|$ and $v$ decreases with respect to $r$.
 \item[iii)]$v\in C^2(\R^d)$.
 \item[iv)] $v$ together with its derivatives up to order 2 have exponential decay at infinity, i.e.
\[
 |D^\alpha v(x)|\leqslant Ce^{-\delta|x|},\quad x\in\R^d,
\]
for some $C, \delta>0$ and for $|\alpha|\leqslant 2$.
\end{itemize}
\end{tm}

Condition \eqref{eqbl2} tells that the existence of a ground state
might depend on the dimension $d$ involved. Indeed, if one
consider the case  $g(v)=-v+|v|^{p-1}v$ and $p>1$, then the
problem \eqref{prob} possesses a solution if
 and only if $1<p<\frac{d+2}{d-2}$. Since we are interested in finding a solution for every $d\geq 3$, we select a different
 nonlinearity, precisely the difference
$|v|^{p_1-1}v - |v|^{p_2-1}v$, with $p_1<p_2$ chosen conveniently,
as explained in the following result (see also \cite[Example
2]{beli}).

\begin{tm}\label{teoex}
 Let $d\geqslant 3$. Then there exist a function $v\in\cS(\R^d)$ and a potential $W(x)\in \cS(\R^d)$ such that
\begin{equation}\label{eqfourth}
 \Delta^2 v - v + W(x)v = 0.
\end{equation}
The potential has the following explicit representation in terms
of the solution $v$
\begin{equation}\label{potA}
 W(x)=-5v^8+24v^6-33v^4+12v^2-20v^2|\nabla v|^2+18|\nabla v|^2.
\end{equation}
\begin{proof}
We shall follow the plan above. Consider the nonlinearity
$g(v)=-v+3v^3-v^5$ (here $p_1=3$ and $p_2=5$). We look for a
solution to
\begin{equation}\label{eqsec}
 -\Delta v + v -3 v^3 + v^5 = 0.
\end{equation}
It's easy to verify that assumptions \eqref{eqbl1},\eqref{eqbl2}
and \eqref{eqbl3} of Theorem \ref{teobl}  are satisfied by $g$.
Observe that with this particular choice of $g$, condition
\eqref{eqbl2} is satisfied in any dimension $d\geqslant 3$. Then,
we come to the  fourth order  PDE by applying the Laplacian
operator to \eqref{eqsec}
\[
\Delta^2 v -\Delta v + 18v|\nabla v|^2+9v^2\Delta v -20v^3|\nabla v|^2 -  5v^4 \Delta v=0.
\]
Finally, we replace $\Delta v$ by its expression in \eqref{eqsec}
and  obtain the desired equation \eqref{eqfourth}, with the
potential $W(x)$ given by \eqref{potA}. The regularity and the
exponential decay of the solution and its derivatives  can be
proved by a bootstrapping argument starting from \eqref{eqfourth}
and using the properties of the solution $v$ in Theorem
\ref{teobl}, Items $iii)$ and $iv)$. Thus we have proved the
existence of a Schwartz solution $v$ to \eqref{eqfourth} that
decay exponentially. Alternatively, the Schwartz and the
exponential decay properties, or, more precisely, the
Gelfand-Shilov class $\cS^{1,1}(\R^d)$ property of the solution
$v$ can be obtained by using \cite[Theorem 4.3]{CGR}.
\end{proof}
\end{tm}

\begin{rem}
Observe that a variational argument like the one developed in
\cite{DaPiVi} is not sufficient to prove Theorem \ref{teoex} for a
given potential $W\in C_0^\infty(\R^d)$. In fact, one would
consider \eqref{eqfourth} as the Euler-Lagrange equation of an
associated minimization problem that doesn't have a minimum by the
following consideration. Let $w(x)$ be a smooth compactly
supported function such that $w(x_0)>0$ at least in one point
$x_0$. Then, if the minimization problem
\[
\min_{f\in M} \int_{\R^d} \big(|\Delta f|^2 - |f|^2 \big) dx \quad \text{on}\quad M=\Big\{f\in H^2: \int_{\R^d}w(x)|f|^2 dx=1\Big\}.
\]
admits a minimum, so the problem \eqref{eqfourth} does. But this
problem has no minimum. It is sufficient to take a function $f_0$
with compact support that satisfies the constraint and a function
$u_1\in H^2$  such that $\int_{\R^d} \big(|\Delta u_1|^2 - |u_1|^2 \big) dx <0$ and having disjoint support from that of  $f_0$. Then
the functions $f_t=f_0+tu_1$, with $t>0$, belongs to $M$ and the
functional $\langle F, f\rangle = \int_{\R^d} \big(|\Delta f|^2 - |f|^2 \big) dx$, applied to $f=f_t$, goes to $-\infty$ as
$t\rightarrow+\infty$.
\end{rem}

\subsection{Sharpness of the global Strichartz estimates}\label{secwell3}

Now we have all the instruments to prove the sharpness of Theorem
\ref{teovib}. We shall detail the global case, the local one can
be obtained by using similar arguments, so we omit the proof.
\begin{tm}
 Let $d\geqslant 3$ and $I=[0,\infty)$. For every $s\in[2,d),\,\alpha\in[1,\infty)$ and  $\beta\in(1,\infty)$, such that
\[
 \frac{2}{\alpha}+\frac{d}{\beta}\neq s+2
\]
there exist a potential $V(t,x)\in L^\alpha \W^{s-2,\beta}$ and a
sequence of solutions $u_k(t,x)\in C_I\Hs^{s}$ to the equation
\eqref{eqomo} such that
\[
\lim_{k\to+\infty}
\frac{\|u_k\|_{L^q\W^{s,r}}}{\|u_k(0)\|_{\Hs^s}+\|\partial_t
u_k(0)\|_{\Hs^{s-2}}}=+\infty, \text{ for every admissible
pair } (q, r)\neq(\infty, 2).
\]
\begin{proof}
The pattern is provided by  \cite[Theorem 1.3]{DaPiVi}. Of course,
working with a higher order equation and with Sobolev spaces
instead of simply Lebesgue spaces, some suitable modifications are
needed. We detail the case
$\displaystyle\frac{2}{\alpha}+\frac{d}{\beta}< s+2$, the case
$\displaystyle\frac{2}{\alpha}+\frac{d}{\beta}> s+2$  is obtained
similarly.

Considering the functions $v$ and $W$ of Theorem \ref{teoex}, we
define the standing wave $u(t,x)=e^{it}v(x)\in
C([0,\infty);L^2(\R^d))$ and observe that, by \eqref{eqfourth},
this function solves the equation \eqref{eqomo} with
\[
 V(t,x)= W(x)\in L^\infty([0,\infty);\W^{s-2,\beta}),
\]
for all $s\in[2,\infty)$ and $\beta\in(1,\infty)$. From Remark
\ref{remresc} we see that the rescaled function
\begin{equation}\label{dil}u_\epsilon(t,x)=e^{i\eps^2t}v(\eps x)
\end{equation}
 solves globally
\[
 \partial_t^2 u_\epsilon + \Delta^2 u_\epsilon + W_\epsilon(x)u_\epsilon = 0, \qquad \text{where}\quad W_\eps(x)=\eps^4 W(\eps x).
\]
Choose two  monotone sequences of positive real numbers
\[
 0=T_0<T_1<\dots<T_k \to +\infty, \qquad \,\,0<e_k\to 0, \qquad k\in\mathbb{N}
\]
and define a potential $V(t,x)$ on $[0,\infty)\times\R^d$  as
follows
\[
 V(t,x)=W_{\eps_k}(x)  \quad \text{for} \quad t\in[T_k,T_{k+1}),\quad k\in\mathbb{N}.
\]
With this choice, $u_{\eps_k}$ solves  the equation
\begin{equation}\label{eqomo2}
 \partial_t^2u+\Delta^2 u+ V(t,x)u=0
\end{equation}
on $[T_k,T_{k+1})$. Select now $\alpha$ and $\beta$ such that
\begin{equation}\label{eqcount}
 \frac{2}{\alpha}+\frac{d}{\beta}<s+2
\end{equation}
and assume we can choose the parameters $T_k$ and $\eps_k$ in such a way that
\begin{equation}\label{eqpar1}
 \|V\|_{L^\alpha \W^{s-2,\beta}}\lesssim \|W\|_{\W^{s-2,\beta}}\sum_{k=0}^\infty\eps_k^{s+2-\frac{d}{\beta}}(T_{k+1}-T_k)^{1/\alpha}<+\infty.
\end{equation}
Then $V\in L^\alpha([0,\infty), \W^{s-2,\beta}(\R^d))$. On the
other hand, Theorem \ref{teovib} allows to extend uniquely the
solutions $u_{\eps_k}$ to a global one of \eqref{eqomo2} in
$C([0,\infty);\Hs^s(\R^d))$ which we shall denote by $u_k(t,x)$.
Observe that we have
$$\| u_k(t,\cdot)\|_{\Hs^s}\equiv \textit{const.}\equiv \|u_{\eps_k}(T_k,\cdot)\|_{\Hs^s}\equiv
\eps_k^{s-\frac d2}\|v\|_{\Hs^s},
$$
where in the last equality we have used the dilation properties
for homogeneous Sobolev spaces. Now, the  Strichartz estimates
\eqref{vistrisolhom} are violated when
\begin{equation}\label{bo}
 \frac{\|u_k\|_{L^q\W^{s,r}}}{\|u_k(0)\|_{\Hs^s}+\|\partial_t u_k(0)\|_{\Hs^{s-2}}} \quad \text{is unbounded}.
\end{equation}
The numerator can be estimated from below by
\[
 \|u_k\|_{L^q\W^{s,r}}\geqslant \eps_k^{s-d/r}\|v\|_{\W^{s,r}}(T_{k+1}-T_k)^{1/q}
\]
while the denominator is estimated from above by
\[
 \|u_k(t,\cdot)\|_{\Hs^s}+\|\partial_t u_k(t,\cdot)\|_{\Hs^{s-2}}=\|v(\eps_k\cdot)\|_{\Hs^s}+\|i\eps_k^2v(\eps_k\cdot)\|_{\Hs^{s-2}}=\eps_k^{s-d/2}(\|v\|_{\Hs^{s}}+\|v\|_{\Hs^{s-2}}).
\]
Thus, the claim \eqref{bo} holds provided that $T_k$, $\eps_k$
satisfy the condition
\begin{equation}\label{eqpar2}
\begin{split}
  \frac{\|u_k\|_{L^q\W^{s,r}}}{\|u_k(0)\|_{\Hs^s}+\|\partial_t u_k(0)\|_{\Hs^{s-2}}} \geqslant & \frac{\eps_k^{s-d/r}\|v\|_{\W^{s,r}}(T_{k+1}-T_k)^{1/q}}{\eps_k^{s-d/2}(\|v\|_{\Hs^s}+\|v\|_{\Hs^{s-2}})}\\
 =& C\eps_k^{d/2-d/r}(T_{k+1}-T_k)^{1/q}\to +\infty
\end{split}
\end{equation}
Finally, by taking
\[
 T_0=0,\quad T_{k+1}=T_k+k^a\quad \text{ and } \quad \eps_0=1,\quad \eps_k=k^{-b/2}, \quad k\in\mathbb{N}
\]
for some $a,b>0$, the conditions \eqref{eqpar1} and \eqref{eqpar2}
reduce to
$$\frac a \alpha+ \frac{b\,d}{2\beta}< \frac{s\, b}2+b-1,\quad  \frac
a q+\frac{b\, d}{2 r}>\frac { b\, d}4.
$$
Since $(q,r)$ is admissible, the second condition becomes $a>b$,
hence, using $a>b$ and the assumptions \eqref{eqcount}, the
first one becomes $\displaystyle\frac{1}{\alpha}+\frac{d}{2\beta}<
\frac s 2+1$.
Thus, the conditions \eqref{eqpar1} and \eqref{eqpar2} are satisfied for
\[
 a>b>2\Big[2-\Big(\frac{2}{\alpha}+\frac{d}{\beta}-s\Big)\Big]^{-1}.
\]
Note that the term in the square bracket is never zero by the
particular assumption on $\alpha,\beta$, see \eqref{eqcount}.
\end{proof}
\end{tm}

\section*{Acknowledgements}
The authors would like to thank Professors Piero D'Ancona, Fabio
Nicola, Luigi Rodino and Paolo Tilli for fruitful conversation and
comments. We are grateful to the anonymous referee for his
valuable comments.

\end{document}